\documentclass[12pt,a4paper,twoside, english]{amsart}

\usepackage[left=2.5cm,right=3.5cm]{geometry}

\usepackage{amssymb}
\usepackage{enumitem}
\usepackage{amsmath}
\usepackage{mathtools}
\usepackage{amsfonts}

\usepackage[utf8]{inputenc}

\usepackage[colorinlistoftodos, textwidth=3cm]{todonotes}

\newtheorem{theorem}{Theorem}[section]

\newtheorem{Lemma }[theorem]{Lemma }

\newtheorem{corollary}[theorem]{Corollary}

\theoremstyle{definition}

\newtheorem{definition}[theorem]{Definition}

\newtheorem{claim}[theorem]{Claim}

\newtheorem{remark}[theorem]{Remark}

\newtheorem*{theorem*}{Theorem}
\newtheorem*{Lemma *}{Lemma }

\newcommand{\ack}{{\bf Acknowledgements}}

\newcommand{\mF}{\mathbb F}
\newcommand{\mE}{\mathbb E}
\newcommand{\mN}{\mathbb N}
\newcommand{\mC}{\mathbb C}

\newcommand{\mT}{\mathbb T}

\newcommand{\mP}{\mathbb P}
\newcommand{\mD}{\mathbb D}

\makeatletter
\newcommand{\dotDelta}{{\vphantom{\Delta}\mathpalette\d@tD@lta\relax}}
\newcommand{\d@tD@lta}[2]{%
  \ooalign{\hidewidth$\m@th#1\mkern-1mu\cdot$\hidewidth\cr$\m@th#1\Delta$\cr}%
}
\makeatother

\usepackage{amsmath}
\usepackage{graphicx}
\usepackage{mathtools}
\usepackage[colorlinks=true, allcolors=blue]{hyperref}
\usepackage{setspace}

\title[Orthogonality of non-pretentious functions to polynomial phases]{Orthogonality of non-pretentious functions to polynomial phases over Function Fields}
\author{Tal Meilin}
\makeatother

\begin{document}

\begin{abstract}
    We prove that non-pretentious multiplicative functions are orthogonal to polynomial phases over $\mF_q[x]$ (up to characteristic conditions). 
\end{abstract}
\maketitle
\tableofcontents
\section{Introduction}
A function $\nu:\mathbb{F}_q[x]\rightarrow\mathbb{C}$ is called multiplicative if $\nu(1)=1$ and for every pair of coprime polynomials $g, h \in \mathbb{F}_q[x]$, it satisfies $\nu(gh)=\nu(g)\nu(h)$. In this paper, we study the correlation of such functions with polynomial phases, a problem deeply connected to the Gowers uniformity of multiplicative functions and the ``pretentious'' framework of Granville and Soundararajan (\cite{MN},\cite{SrGr}).

The problem of multiplicative functions correlating with polynomial phases has its roots in the classical setting of the integers $\mathbb{Z}$. An important result in this area is the work of Green and Tao \cite{MN}, who proved that the M\"{o}bius function is orthogonal to Nilsequences. Green and Tao's result provided an effective bound of $O(\log^{-A} N)$, which is the limit of unconditional results in $\mathbb{Z}$ due to potential Siegel zeros. If one assumes the Generalized Riemann Hypothesis (GRH), this bound improves to a power-saving $O(N^{-c})$.

In the context of function fields, we proved the analogous result to Green and Tao's work for the M\"{o}bius function, under certain characteristic conditions (\cite{mobius}).

The effectiveness of the bound in \cite{mobius} relies on the use of Vaughan's Lemma, which is a structural identity specific to the M\"{o}bius function.

The goal of this paper is to generalize these results to a much larger class of multiplicative functions --- specifically, those that are aperiodic.
\begin{definition}[Hayes-Periodic bounded Function]
    A bounded function $a:\mF_q[x]\rightarrow\mC$ is Hayes-periodic if there exist some $Q\in\mF_q[x]$ and $L\in\mN$ such that for every $f,g\in\mF_q[x] $ w:
    \begin{enumerate}
        \item $f\equiv g \text{ mod (Q)}$
        \item $f,g$ has the same $L$ highest coefficients
    \end{enumerate}
    we have $a(f)=a(g)$.
\end{definition}
\begin{definition}[Aperiodic Function] A bounded multiplicative function $\nu: \mathbb{F}_q[x] \to \mathbb{C}$ is said to be aperiodic if it is orthogonal to every Hayes-periodic function $a:\mF_q[x]\rightarrow\mC$ in the sense that:$$\mathbb{E}_{g,\deg(g)=n} \nu(g) \overline{a}(g) = o(1) \quad \text{as } n \to \infty.$$
 
\end{definition}
\begin{remark}
    Observe that in particular, if $\nu$ is aperiodic then:
    \[\mathbb{E}_{g \text{ monic},\deg(g)=n} \nu(g) \overline{a}(g) = o(1) \quad \text{as } n \to \infty\]
\end{remark}
Our main result establishes that aperiodicity alone is sufficient for orthogonality to polynomial phases of degree $<p$. 

\begin{theorem}[Main Theorem]\label{main}
 Let $q=p^r$, and let $\nu: \mathbb{F}_q[x] \to \mathbb{C}$ be a $1$-bounded aperiodic multiplicative function. Then for any additive polynomial $P: \mathbb{F}_q[x] \to \mathbb{F}_q$ of degree $d < \operatorname{char}(\mathbb{F}_q)$ (see Definition \ref{poly}), and any additive character $\alpha$:
$$\sum_{g \in G_n} \nu(g) \alpha(P(g)) = o_{\nu,d}(q^n),$$   
and the bound is uniform in $P$.
\end{theorem}

\begin{remark}
    As mentioned in Lemma \ref{aperiodic-characters}, a function is aperiodic if and only if it is uniformly orthogonal to the exponential of every linear (additive) polynomial. Theorem \ref{main} extends this criterion to orthogonality to every polynomial phase of degree $< \operatorname{char}(\mathbb{F}_q)$.
\end{remark}
While this result is ineffective (providing a $o(1)$ type bound, rather than $q^{-n^c}$), it applies to {\em all} aperiodic, bounded multiplicative functions.

Note that while the result in \cite{mobius} was under the condition $d<\frac{p}{2}$, it can be improved to $d<p$ by using the same methods we used in this paper.

Similarly to the context of $\mN$, this result is intimately related to the proximity of an aperiodic function $\nu$ to certain multiplicative characters. 

We define the pretentious distance between two multiplicative functions $f, g$ as:
$$\mathbb{D}(f, g, n) = \left( \sum_{r\in\mP,\ \deg(r) \le n} \frac{1 - \text{Re}(f(r)\overline{g}(r))}{q^{\deg(r)}} \right)^{1/2}.$$

In the context of the integers $\mN$, a function is aperiodic if and only if its pretentious distance from all twisted Dirichlet characters (i.e., functions of the form $\chi(n)n^{it}$ where $\chi$ is a Dirichlet character) is infinite.

In the context of $\mF_q[x]$ we will use a similar criterion, but instead of twisted Dirichlet characters, we will measure the pretentious distance from Hayes characters (Definition \ref{hc}). Similar ideas were presented in the work of Klurman, Mangerel, and Ter\"av\"ainen in \cite{PDT}, as well as the function fields version of the Halasz-Delange Theorem (Theorem \ref{haldel}).

In fact, Theorem \ref{main} follows from the following more general theorem:

\begin{theorem}\label{eq}
For a $1$-bounded multiplicative function $\nu: \mathbb{F}_q[x] \to \mathbb{C}$, the following are equivalent:
\begin{enumerate}
    \item $\nu$ is aperiodic.
    \item $\nu$ is uniformly orthogonal to every polynomial phase $\alpha(P)$ where $P$ is an additive polynomial of degree $d < \operatorname{char}(\mathbb{F}_q)$ and $\alpha$ is an additive character.
    \item For each $K,l$ \[\lim_{n \to \infty} \min_{\xi,cond(\xi)\leq l}\min_{\chi,cond(\chi)\leq K}\min_{\theta\in[0,1)} \mathbb{D}(\nu, \xi\chi e_\theta, n) = \infty,\] 
    where $\chi$ is a Dirichlet characters, $\xi$ is a short interval characters.
\end{enumerate}
\end{theorem}

By invoking the Inverse Theorem of the Gowers Norms (\cite{IGN}), Theorem \ref{eq} can be viewed as the function fields analogue of the result of Frantzikinakis and Host:
\begin{theorem}[\cite{HOFA}, Theorem 2.5]\label{Fr-Ho}
    If a bounded multiplicative function $f$ is aperiodic, then $\lim_{N\rightarrow\infty}\|f\|_{U^s[N]}=0$ for every $s\geq2$.
\end{theorem}
Combining the Inverse Theorem of the Gowers Norms with a result by Gowers (\cite{Gowers}, Theorem 3.2), we deduce the following:
\begin{theorem}\label{AP}
     Let $p$ be prime and $\nu: \mathbb{F}_p[x] \to \mathbb{C}$ be a $1$-bounded, aperiodic multiplicative function. Then for every $k < p$ and $1$-bounded functions $\nu_1, \dots, \nu_k$:$$\left|\mathbb{E}_{x, y \in G_n} \nu_1(x) \nu_2(x+y) \cdots \nu_k(x+(k-1)y) \nu(x+ky)\right| = o(1).$$
     In particular, we may choose $\nu_1=\nu_2=\ldots=\nu_{k}=\nu$ and get that:
    \[|\mE_{x, y \in G_n}\nu(x)\nu(x+y)\cdot\ldots\cdot \nu(x+ky)|=o(1).\]
\end{theorem}

\ack. I thank Joni Ter\"av\"ainen for his helpful ideas and suggestions. I am especially grateful to my academic advisor Tamar Ziegler, for her guidance and support throughout this project. This work was supported by ISF grant 2112/20.

We used AI tools for language editing and to assist with checking proofs and performing some calculations related to Theorems \ref{aperiodic-characters} and \ref{apf-equiv-np}.

\section{Notation}
Let $q=p^r$ be the power of a prime number $p$.
We begin by defining the set of additive characters of $\mathbb{F}_q$, which plays an important role in Theorem \ref{main}.

Recall that each additive character of $\mathbb{F}_q$ is of the form:
\[t\xrightarrow{}\exp\left(\frac{2 \pi i \cdot \operatorname{Tr}(t s)}{p}\right),\]
for some $s\in \mathbb{F}_q$, where $\operatorname{Tr}:\mathbb{F}_q\xrightarrow{} \mathbb{F}_p$ is the trace map.

Define the field of formal Laurent series as:$$\mathbb{F}_q((1/x))=\left\{\alpha=\sum_{i=-\infty}^{n} a_i x^i \mid n\in \mathbb{Z},  a_i\in\mathbb{F}_q \text{ for every } i\right\}.$$

We define a norm on the field $\mathbb{F}_q((1/x))$ by $\|\alpha\|=q^{\deg (\alpha)}$. Note that $\deg(\alpha)$ might be negative. Define the function fields analogue of the unit interval as:$$\mathbb{T}=\left\{\alpha\ \in \mathbb{F}_q((1/x))\mid \|\alpha\|<1\right\}.$$

If $\mathbb{F}_q[x]$ is analogous to $\mathbb{Z}$, then $\mathbb{F}_q((1/x))$ is analogous to $\mathbb{R}$, and $\mathbb{T}$ is analogous to $\mathbb{R} / \mathbb{Z}$. In this setting, the analogue to $\mathbb{Q}$ is $\mathbb{F}_q(x)$ --- the field of fractions of $\mathbb{F}_q[x]$. 

Similarly to $\mathbb{R}$, we say that an element $a\in\mathbb{F}_q((1/x))$ is rational if $a\in\mathbb{F}_q(x)$ and irrational if $a\notin\mathbb{F}_q(x)$. 

For each $f\in \mathbb{F}_q((1/x))$, we denote by $f_{-1}$ the coefficient of $x^{-1}$. Note that each linear form of $\mathbb{F}_q[x]$ is of the form $f\xrightarrow{} (\beta f)_{-1}$ for some $\beta\in\mathbb{T}$.

\subsection{Derivatives, Additive Polynomials and Multilinear maps}
Throughout this paper, we will use the following definitions:
\begin{definition}[Additive Polynomials]\label{poly} Let $V$ be a $\mathbb{F}_q$-vector space. A function $P:V\rightarrow\mathbb{F}_q$ is an additive polynomial of degree $\leq m$ if for every $h_0,\ldots,h_m\in V$, $P$ satisfies:$$\Delta_{h_0}\ldots\Delta_{h_m}P\equiv 0,$$where $\Delta_hP(g)=P(g+h)-P(g)$.
\end{definition}

We will now turn to the relation between additive polynomials and multilinear maps.

\begin{definition} Let $Q:\mathbb{F}_q[x]^m \to \mathbb{F}_q$ be an $m$-$\mF_p$linear map. Define $P_Q$ to be the homogeneous additive polynomial corresponding to $Q$, obtained by restricting $Q$ to the diagonal:$$P_Q(g):=Q(g,\ldots,g),\text{ for every } g \in \mathbb{F}_q[x].$$\end{definition}

Furthermore, we define the multilinearization $d^m H$ of an additive polynomial:

\begin{definition}Let $H:\mathbb{F}_q[x]\rightarrow\mathbb{F}_q$ be an additive polynomial of degree $\leq m$. Define $d^mH:\mathbb{F}_q[x]^m\rightarrow\mathbb{F}_q$ by:$$d^mH(h_1,\ldots,h_m):=\Delta_{h_1}\ldots\Delta_{h_m}H(g).$$It is straightforward to verify that $d^mH$ is an $m$-$\mF_p$ linear map independent of the choice of $g$.
\end{definition}

 Another important definition for our paper is phase polynomials:
 \begin{definition}
     Let $\alpha\in\hat{\mF_q}$, and let $\tilde{H}:\mathbb{F}_q[x]\rightarrow\mathbb{F}_q$ be an additive polynomial of degree $\leq m$. Then we say that $\tilde{H}$ is a phase polynomial of degree $\leq m$. In particular, for convenience, we will use extensively the specific character:
     \[e_p(s)=exp(\frac{2 \pi i}{p}Tr_{\mF_q/\mF_p}(s)).\]
     Under this convention, will also considered the polynomial $H=Tr(s\tilde{H})$. 
 \end{definition}

The relation between additive polynomials and multilinear maps is characterized by the following identities:
\begin{claim}Let $H:\mathbb{F}_q[x]\rightarrow\mathbb{F}_q$ be an additive homogeneous polynomial of degree $ m$, and $Q:\mathbb{F}_q[x]^m \to \mathbb{F}_q$ be an $m$-$\mF_p$ linear map. Then:
\begin{enumerate}
    \item $P_{d^m H}=m! H$.
    \item $d^mP_Q=m! Q$.
\end{enumerate}\end{claim}

\subsection{Schmidt Rank and Partition Rank}The concept of rank is used to determine how easily a particular polynomial can be expressed as a function of lower-degree polynomials. Since Theorem $\ref{main}$ is proven through induction, expressing high-degree polynomials as functions of simpler ones is essential.

\begin{definition}[Schmidt Rank]\label{AR}Let $V$ be a $\mF_q$-vector space. For an additive, homogeneous polynomial $P:V\rightarrow\mF_q$ of degree $m>1$, the Schmidt rank $r(P)$ is the minimum integer $r$ s.t:$$P= \sum_{i=1}^r M_iR_i$$where $M_i$ and $R_i$ are additive, homogeneous polynomials of positive degrees $<m$. For a non-homogeneous polynomial $P$, the Schmidt rank is defined as the rank of its highest-degree homogeneous part.
\end{definition}
Similarly, we define a rank for multilinear forms:
\begin{definition}[Partition Rank]\label{PR}For an $m$-linear map $Q$, the partition rank $r_{pr}(Q)$ is the minimal $r$ s.t:$$Q= \sum_{i=1}^r M_iR_i$$where $M_i, R_i$ are $k_i$-linear and $(m-k_i)$-linear maps respectively, with $1 \le k_i < m$.\end{definition}

\begin{remark} \label{r-ineqaulity}
Let $Q$ be a $m$-linear map, and $P$ be an additive polynomial of degree $m$. Then:
\begin{enumerate}
    \item If $m<\operatorname{char}(\mF_q)$ and $r_{pr}(Q)=s$, then $r(P_Q)\leq s$.
    \item If $r(P)=r$ then $r_{pr}(d^mP)\leq2^m r(P)$.
\end{enumerate}
\end{remark}
\begin{proof}
    The first statement is trivial, based on the definition of rank and the definition of $P_Q$.

For the second one, assume $r(P)=r$:    
\[d^mP(h_1,\ldots, h_m)=\sum_{i=1}^r \Delta_{h_1}\ldots \Delta_{h_m} M_i R_i,\]
where $M_i\text{ and }R_i$ are homogeneous of degree $k_i, m-k_i$ for some $1\leq k_i<m$ respectively. Hence:
\begin{align*}
    \Delta_{h_1}\ldots \Delta_{h_m} M_i R_i=\sum_{I\subset[m], |I|=k_i}\Delta_I M_i \Delta_{I^c} R_i,
\end{align*}
where:
\[\Delta_{I}M=\Delta_{h_{i_1}}\ldots\Delta_{h_{i_k}}M,\]
where $i_1,\ldots, i_k$ are the elements of $I.$ 
Thus:
\[d^mP(h_1,\ldots, h_m)=\sum_{i=1}^r \Delta_{h_1}\ldots \Delta_{h_m} M_i R_i=\sum_{i=1}^r\sum_{I\subset[m], |I|=k_i}\Delta_I M_i \Delta_{I^c} R_i, \]
hence:
\[r_{pr}(d^mP)\leq2^m r(P).\]
\end{proof}

\begin{remark}
    In this paper, we often start with an additive polynomial $$P:\mF_q[x]\rightarrow\mF_q,$$ and identify $P$ with its restriction to a large enough finite-dimensional subspace of $\mF_q[x]$ (for example, $\tilde{V}=\{g\in\mF_q[x]\mid \deg(g)\leq n\}$ for large enough $n$).
\end{remark}

\subsection{Multiplicative Characters and Pretentious Distance}In the first part of the paper, we examine the correlation of multiplicative functions with phase polynomials. The second part turns to multiplicative ``obstructions'', specifically Hayes characters.

\begin{definition}[Dirichlet characters]\label{dc}
A completely multiplicative function $\chi:\mathbb{F}_q[x]\rightarrow\mathbb{C}$ is a {\em Dirichlet character} if there exists $g\in\mathbb{F}_q[x]$ such that for every $h\in\mF_q[x]$, $\chi$ satisfies $\chi(h)=\chi(h \bmod g)$ and $\chi(h) \neq 0$ if and only if $\gcd(g,h)=1$.
\end{definition}

Over the integers, twisted Dirichlet characters take the form $\chi(n)n^{it}$. In $\mathbb{F}_q[x]$, the role of $n^{it}$ is played by degree twists $e_{\theta}(g) = \exp(2i\pi \theta \deg(g))$, but Dirichlet characters alone are not enough. We also need to consider short interval characters.

\begin{definition}[Short interval characters]\label{sic}
A multiplicative function $\xi:\mathbb{F}_q[x]\rightarrow\mathbb{C}$ is a {\em short interval character} of length $s$ if $\xi(g)$ depends only on the $s+1$ highest degree coefficients of $g$.
\end{definition}

\begin{definition}[Hayes characters]
\label{hc}A {\em Hayes character} $\tilde{\chi}=\chi\xi$ is a product of a Dirichlet character $\chi$ and a short interval character $\xi$.
\end{definition}

\begin{definition}[Conductor]
    For a Dirichlet character $\chi$, the {\em conductor} of $\chi$, $cond(\chi)$, is the degree of an element $g\in\mF_q[x]$ such that:
    \begin{enumerate}
        \item $\chi(h)=\chi(h \bmod g)$ and $\chi(h) \neq 0$ if and only if $\gcd(g,h)=1$.
        \item for any $h'$ with the previous property, $h|h'$. 
    \end{enumerate}
    Similarly, for a short interval character $\xi$, $cond(\xi)$ is the smallest $s$ such that for every $\xi(g)$ depends only on the $s+1$ highest degree coefficients of $g.$\\
    Finally, for a Hayes character $\tilde{\chi}=\chi\xi$:
    \[cond(\tilde{\chi})=cond(\chi)+cond(\xi).\]
\end{definition}
To measure how much a multiplicative function ``pretends'' to be a Hayes character, we use the following definition of distance:
\begin{definition}[Pretentious distance]\label{pd}
Let $f,g :\mathbb{F}_q[x]\rightarrow\mathbb{C}$ be bounded multiplicative functions. We define the {\em pretentious distance} between $f,g$ as:$$\mathbb{D}(f,g,N)=\left(\sum_{p\in\mathbb{P}_{\leq N}}q^{-\deg(p)}\left(1-\operatorname{Re}(f(p)\bar{g}(p))\right)\right)^{0.5}.$$
\end{definition}

Finally, we characterize non-pretentious functions:
\begin{definition}[Non-pretentious function]
A $1$-bounded multiplicative function $f$ is {\em non-pretentious} if for every $K$:
\[\lim_N \min_{\tilde{\chi},cond(\chi)\leq K}\min_{\theta\in[0,1]} \mathbb{D}(f,\tilde{\chi} e_{\theta},N)=\infty,\]
where the minimum of $\tilde{\chi}$ is taken over the Hayes-characters of conductor $\leq K$.
\end{definition}

\section{Preliminaries}
As mentioned, we prove Theorem \ref{main} through induction on the degree $k$ of the polynomial phase. Assuming by contradiction that the theorem fails for some $k < \operatorname{char}(\mathbb{F}_q)$ and a polynomial phase $P$ of degree $k$, we aim to derive structural properties about $P$ that lead to a contradiction. To achieve this, we rely on the following key tools.

The first is a function-field adaptation of a classical result in the integer setting. Katai's criterion lets us reduce bounding $\sum \nu(g)f(g)$ to controlling sums $\sum f(ag)\overline{f}(bg)$ over pairs of irreducible polynomials of the same degree. In our setting with $f=\alpha\circ P$, the resulting condition on $f$ is what forces $P$ to have bounded Schmidt rank.

\begin{theorem}[Katai's criterion, proved in Section~\ref{appendix}]\label{p_katai}
Let $f:\mF_q[x]\to\mathbb{C}$ with $|f|\le1$, and let $\nu:\mF_q[x]\to\mathbb{C}$ be
multiplicative with $|\nu|\le1$. Let $0\le W<H$ be integers with $2H\le n$, and let $\eta\ge0$
be such that for every integer $k$ with $W<k<H$,
\begin{equation}\label{eq:katai-hyp}
\sum_{a,b\in P_k}\Big|\sum_{g\in G_{n-k}}f(ag)\overline{f}(bg)\Big|\le \eta\,|P_k|^2q^{n-k}.
\end{equation}
Then
\[
\Big|\sum_{g\in G_n}\nu(g)f(g)\Big|\le q^n\Big(A^{-1/2}+\frac{2q^{-W}}{(W+1)A}+\eta^{1/2}\Big).
\]
\end{theorem}

Theorem \ref{p_katai} allows us to convert a correlation problem involving a multiplicative function into a structural problem concerning the polynomial phase $f=\alpha(P)$. Specifically, we use this to show that significant correlation with an aperiodic function implies that $P$ has low rank. 
Consequently, we require the following result regarding the behavior of Schmidt rank under restriction to a subspace, as found in \cite{k-z}.

\begin{Lemma }[\cite{k-z}, Remark 3.7]\label{sub} Let $\mathbb{F}$ be a field, $V$ an $\mathbb{F}$-vector space, and $P: V \to \mathbb{F}$ an additive polynomial of degree $m$. Let $W \subset V$ be a subspace of codimension $k$. If $P_W$ denotes the restriction of $P$ to $W$, then:$$ r(P_W) \geq r(P) - mk. $$
\end{Lemma }

The primary analytical tool for proving low-rank properties from exponential sums is the Bias-rank theorem. It establishes that if an additive polynomial (or multilinear map) is biased (meaning, its exponential sum does not vanish), it must be expressible via a small number of lower-degree terms.

\begin{theorem}[Bias-rank, \cite{BRc} Theorem 3] \label{BR}
Let $V$ be an $n$-dimensional vector space over $\mathbb{F}_q$, $Q$ an $m$-linear map on $V^m$, and $\alpha_1: \mathbb{F}_q \to \mathbb{C}$ the additive character $\alpha_1(a) = \exp(2\pi i \operatorname{Tr}(a)/p)$. Then:
\begin{enumerate}
    \item $\mathbb{E}_{x \in V^m} \alpha_1(Q(x)) \geq q^{-r_{pr}(Q)}$. In particular, the expectation is always non-negative.
    \item If $\mathbb{E}_{x \in V^m} \alpha_1(Q(x)) \geq q^{-s}$, then $r_{pr}(Q) < C(m)s^{C(m)}$, where $C(m)$ depends only on $m$.
\end{enumerate}
This conclusion holds similarly for an additive polynomials $P$ of degree $m < \operatorname{char}(\mathbb{F}_q)$: if $$|\mathbb{E}_{x \in V} \alpha_1(P(x))| \geq q^{-s},$$
then $r(P) < C(m)s^{C(m)}$.
\end{theorem}



\section{Proof of the main Theorem}

In this section, we provide the proof for Theorem \ref{main}. We proceed by induction on the degree $m$ of the polynomial phase.

The starting point of our induction is the relationship between aperiodicity and linear phases. We observe that the property of aperiodicity, originally defined by orthogonality to periodic sequences, is equivalent to cancellation against linear polynomial phases:

We start by proving that aperiodic functions are uniformly orthogonal to linear maps, i.e to the characters of $\mF_q[x]$.
\begin{theorem}\label{aperiodic-characters}
    Let $\nu:\mF_q[x]\rightarrow\mC$ a bounded multiplicative function. Then $\nu$ is aperiodic if and only if:
    \[\sup_{\theta\in\mT}|\mE_{f\in G_N}\nu(g)e_{\mF}(\theta g)|=o(1),\]
    where $e_{\mF}(\theta g)=exp(2\pi i\cdot Tr_\mF((\theta g)_{-1}))$.
\end{theorem}
\begin{proof}
\textit{ Apriodicity $\implies$ Uniform Orthogonality}
Observe that:
\begin{align*}
    &\left|\sum_{f\in G_N}\nu(g)e_{\mF}(\theta g)\right|=\left|\sum_{n=0}^{N-1}\sum_{g,\deg(g)=n}\nu(g)e_{\mF}(r\theta g)\right|
    &\leq \sum_{n=0}^{N-1}\left|\sum_{g,\deg(g)=n}\nu(g)e_{\mF}(r\theta g)\right|
\end{align*}
Hence it is enough to prove that for $n$ large enough:
\[\sup_{\theta\in\mT}\left|\sum_{g,\deg(g)=n}\nu(g)e_{\mF}(r\theta g)\right|=o(q^n).\]
Define a new norm on $\mF_q((\frac{1}{x}))$:
\[\|\theta\|_{frac}=q^{\deg(\theta-\lfloor \theta\rfloor)}.\] We will use this definition in this chapter.
    Fix some $R,L>0$, and assume $n>2R,2L$. Split $\mT$ into two disjoint set:
    \begin{itemize}
        \item Major Arcs: $\theta$ is a major arcs if there exists some $b\in\mF_q[t]$ with $\deg(b)<R$ and $\|\theta b\|_{frac}\leq q^{L-n-1}$. \\
        In other word, we can find $a,b\in\mF_q[t]$ with $\deg(a),\deg(b)<R$ and $\beta\in\mT$ with $\deg(\beta)\leq L-n-1$ such that:
        \[\theta=\frac{a}{b}+\beta\]
        \item Minor Arcs: $\theta$ is a minor arc if for every $b\in\mF_q[t]$ with $\deg(b)<R$, $\|\theta b\|_{frac}>q^{L-n-1}$.
    \end{itemize}

\emph{Case 1: $\theta$ is a major arc:}\\
Observe that for $f,g$ monics of degree $n$, if $f\equiv g \text{ mod b}$ and $f,g$ has the same top $L+1$ coefficients, then:
\[e_{\mF}(\theta g)=e_{\mF}(\frac{a}{b}g)e_{\mF}(\beta g).\]
Since $f\equiv g \text{ mod b}$ we have:
\begin{align*}
    &e_{\mF}(\frac{a}{b}g)=exp\left(\frac{2\pi i}{p}Tr\left((\frac{a}{b}g)_{-1}\right)\right)=exp\left(\frac{2\pi i}{p}Tr\left((\frac{a}{b}f)_{-1}\right)\right),
\end{align*}
and since $f,g$ has the same $L+1-$largest coefficients and $\deg(\beta)<L-n-1$ we have:
\begin{align*}
    &e_{\mF}(\beta g)=exp\left(\frac{2\pi i}{p}Tr\left((\beta g)_{-1}\right)\right)=exp\left(\frac{2\pi i}{p}Tr\left(\sum_{i=0}^{n}g_i\beta_{-i-1}\right)\right)\\
    &=exp\left(\frac{2\pi i}{p}Tr\left(\sum_{i=n-L}^{n}g_i\beta_{-i-1}\right)\right)\\
    &=exp\left(\frac{2\pi i}{p}Tr\left(\sum_{i=n-L}^{n}f_i\beta_{-i-1}\right)\right)=exp\left(\frac{2\pi i}{p}Tr\left((\beta f)_{-1}\right)\right)=e_{\mF}(\beta f)
\end{align*}
where $g=\sum_{i=0}^ng_ix^i$. Meaning $e_{\mF}(\theta g)$ depends only on $g \text{ mod b}$ and on the top $L-$ coefficients of g, $g_{n-L},\ldots, g_{n}$.

Hence, if we define $[g]_L=\sum_{i=0}^{L}g_{n-L+i}x^{i}$ then:
\begin{align*}
    &\sup_{\theta, \theta \text{ is major arc}}|\mE_{\deg(f)=n}\nu(g)e_{\mF}(\theta g)|\\
    &\leq \max_{b,\deg(b)<R}\max_{A\in G_n/(b)}\max_{h\in G_L} |\mE_{\deg(f)=n, f\equiv A \text{ mod b}, [g]_L=h}\nu(g)|=o(1).
\end{align*}
Since the maximum runs over a finite collection, and $\nu$ is aperiodic.\\
\\
\emph{Case 2: $\theta$ is a Minor Arc:}\\
Apply Theorem~\ref{p_katai} with $f(g)=e(\theta g)$. For $W<k\leq \min\{L-1,R\}$ and $a,b\in \mP_k$ we have
$f(ag)\overline{f}(bg)=e\big((a-b)\theta g\big)$.

If $a\neq b$, then $c=a-b\neq0$ and $\deg c<k\le\min\{ L-1,R\}$, so
$\|c\theta\|_{frac}>q^{L-n-1}\ge q^{-(n-k)}$. 

Hence
$\sum_{g\in G_{n-k}}e(c\theta g)$ vanishes.

Hence only the diagonal contributes:
\[
\sum_{a,b\in P_k}\Big|\sum_{g\in G_{n-k}}f(ag)\overline{f}(bg)\Big|=|\mP_k|\,q^{n-k}
=\frac{1}{|\mP_k|}\,|\mP_k|^2q^{n-k}.
\]
So \eqref{eq:katai-hyp} holds with \[\eta=\max_{W<k<\min\{L-1,R\}}|\mP_k|^{-1}\le\max_{k>W}2k\,q^{-k}=2(W+1)q^{-(W+1)},\]
Theorem~\ref{p_katai} imply that:
\[
\Big|\mE_{g\in G_n}\nu(g)f(g)\Big|\le A^{-1/2}+\frac{2q^{-W}}{(W+1)A}+\eta^{1/2}, \qquad A=\sum_{W<\deg(r)<\min\{R,L-1\}}q^{-\deg(r)}.
\]

Let $\epsilon>0$. Choosing $R,L,W\leq\log(L),\log(R)$ large enough (depending only on $\epsilon$) such that:
\[A^{-1/2}+\frac{2q^{-W}}{(W+1)A}+\eta^{1/2}\leq \epsilon/4,\]
and since $\Big|\sum_{\deg(g)=n}\nu(g)f(g)\Big|=\Big|\sum_{g\in G_{n+1}}\nu(g)f(g)-\sum_{g\in G_n}\nu(g)f(g)\Big|$,
and by choosing $n$ large enough, we are done.\\
\textit{ Apriodicity $\impliedby$ Uniform Orthogonality}: Assume that:
 \[\sup_{\theta\in\mT}|\mE_{f\in G_N}\nu(g)e_{\mF}(\theta g)|=o(1).\] 
 Let $Q\in\mF_q[x], \deg(Q)=K$ and $L\in\mN$, and a Hayes-periodic function $a(g)=a(g \text{ mod }Q, g_{\deg(g)-L},\ldots,g_{\deg(g)})$. 

By definition, there exists $a_{s,h}\in\mC$ such that for every $g$:
\[a(g)=\sum_{s\in \mF_q^L}\sum_{h\in \mF_q[x]/Q}a_{s,h}1_{g\equiv_Qh, (g_{\deg(g)-L},\ldots,g_{\deg(g)})=s}.\]
Hence, it's enough to prove that $\nu$ is orthogonal to every $1_{g\equiv_Qh, (g_{\deg(g)-L},\ldots,g_{\deg(g)})=s}$.

Let $n>K,L$.
Observe that the following maps are linear:
 \[\chi:g\rightarrow g \text{ mod }Q, \xi:g\rightarrow \sum_{i=0}^L g_{n-L+i}t^i .\]
In particular $e_p(\chi),e_p(\xi)$ are characters.
For each $M$ define:
\[D_M=\{\sum_{i=1}^M a_ix^{-i}|\ a_i\in\mF_q\}.\]
Using Fourier transform:
\begin{align*}
    &1_{g\equiv_Qh, (g_{\deg(g)-L},\ldots,g_{\deg(g)})=s}\\
    &=\mE_{\theta\in D_K}e_p(\theta(\chi(g)-h))\mE_{\alpha\in D_{L+1}}e_p(\alpha (\xi(g)-s))\\
    &=\mE_{\theta,\in D_K,\alpha\in D_{L+1}}e_p(\theta\chi(g)+\alpha\xi(g))e_p(-\theta s-\alpha h)
\end{align*}
Meaning, it's enough to show that $\nu$ is orthogonal to  $e_p(\theta(\chi(g)+\alpha\xi(g))e_p(-\theta s-\alpha h)$ on the set of polynomials of degree $n$, for any $\theta,\alpha,h,s$. Since this map is an additive character, we are done.
\end{proof}
Following this characterization of Theorem \ref{aperiodic-characters}, we assume by induction that $\nu$ is uniformly orthogonal to every polynomial phase of degree at most $m-1$, and we aim to prove the result for polynomial phases of degree $m < \operatorname{char}(\mathbb{F}_q)$.

To quantify the inductive hypothesis, we define the maximum correlation for lower-degree additive polynomials.

To prove Theorem \ref{main}, we assume by contradiction that $\nu$ is not orthogonal to some exponent of some additive polynomial $P$ of degree $m$. We aim to show that this non-orthogonality forces $P$ to have a low rank, allowing it to be approximated by lower-degree additive polynomials.\\

\textit{ Proof of theorem \ref{main}:} The following lemma establishes that if $\nu$ correlates with $P$, then a certain multilinear derivative $R$ must exhibit significant bias.
We will need the following definitions:
\begin{definition}
    Define: 
    \[\phi_{m-1}(n)=\sup_{H,deg(H)\leq m-1}|\sum_{f\in G_n}\nu(f)\alpha(H(f))|,\]
 In addition, define:
\[
\Lambda(n)=\min\Big\{\log_q\frac{q^n}{\phi_{m-1}(n)},\ n\Big\},\]  and
\[\
W(n)=\Lambda(n)^{1/8},\qquad H(n)=\Lambda(n)^{1/4},\qquad
A(n)=\sum_{W<\deg r<H}q^{-\deg r}.
\]
\end{definition} 
    
To prove this theorem, we start by proving the following lemma: 
\begin{Lemma }[correlation forces low rank, uniformly]\label{key}
Assume $\phi_{m-1}(n)=o(q^n)$. For every $\delta\in(0,1)$ there are $C_1=C_1(\delta,m,q)$ and
$n_1=n_1(\delta,m,q,\Lambda(\cdot))$ with the following property.

For every $n\ge n_1$, every polynomial $\tilde{P}$ of degree $\le m$ and every nontrivial $\alpha=e_p(t\cdot)$: if
\[
\Big|\mathbb{E}_{g\in G_n}\nu(g)\,\alpha(\tilde{P}(g))\Big|>\delta,
\]
then
\[
r_{pr}(d^mTr(t\tilde{P}))\ \le\ C_1\,(\log\Lambda(n))^{C_{2m}}+mr\,H(n).
\]
In particular, the bound on the rank does not depend on $\tilde{P}$.
\end{Lemma }
\begin{proof}
    Let $\tilde{P},\alpha=exp(\frac{2\pi i}{p}Tr(t\cdot))$ satisfies the conditions of the lemma. Define:
    \[P(g):\mF_q[x]\rightarrow\mF_p,\ P(g)=Tr(t\tilde{P}(g))_{\mF_q/\mF_p}=Tr(t\tilde{P}(g)).\]
    Then $e(P)\equiv\alpha(\tilde{P}),$ where $e(r)=exp(\frac{2\pi i}{p}r).$
    
      Using the induction assumption, we have that $\phi_{m-1}(n)=o(q^n)$, meaning $\Lambda(n),W(n),H(n)$ are all tending to infinity.\\
      Let $n$ be so large that $A^{-1/2}+\frac{2q^{-W}}{(W+1)A}\le\frac\delta2$.
If \eqref{eq:katai-hyp} held for all $k\in(W,H)$ with $\eta=\delta^2/4$, then Theorem \ref{p_katai} would give
$|\mathbb{E}_{G_n}\nu\,e(P)|\le\delta$. Hence there is $k\in(W,H)$ with
\[
\mathbb{E}_{a,b\in P_k}\Big|\mathbb{E}_{g\in G_{n-k}}e\big(P(ag)-P(bg)\big)\Big|>\frac{\delta^2}{4}.
\]
(Here $ag\in G_n$, since $\deg(ag)\le k+(n-k-1)=n-1$.) Consequently the set
\[
\tilde \mP=\Big\{(a,b)\in \mP_k^2:\ \big|\mathbb{E}_{g\in G_{n-k}}e(P(ag)-P(bg))\big|\ge\tfrac{\delta^2}{8}\Big\}
\]
satisfies $|\tilde \mP|\ge\frac{\delta^2}{8}|\mP_k|^2$.\\

 For $a,b\in G_{k+1}$ and $h\in G_{n-k}^m$ put
\[
P_{a,b}(g)=P(ag)-P(bg),\qquad
R(a,b,h)=d^mP(ah_1,\dots,ah_m)-d^mP(bh_1,\dots,bh_m).
\]
Since $g\mapsto ag$ is linear $G_{n-k}\to G_n$, $d^m(P_{a,b})=R(a,b,\cdot)$, and hence
\[
\mathbb{E}_{h\in G_{n-k}^m}e(R(a,b,h))=\|e\circ P_{a,b}\|_{U^m(G_{n-k})}^{2^m}\ \ge0 .
\]
For $(a,b)\in\tilde \mP$, monotonicity of Gowers norms gives
\[\|e\circ P_{a,b}\|_{U^m}^{2^m}\ge\|e\circ P_{a,b}\|_{U^1}^{2^m}\ge(\delta^2/8)^{2^m}.\]
All terms are nonnegative, and $|\mP_k|\ge q^k/(2k)$ with $k<H$. Keeping only $\tilde P$,
\[
\beta_1:=\mathbb{E}_{a,b\in G_{k+1}}\ \mathbb{E}_{h\in G_{n-k}^m}e(R(a,b,h))
\ \ge\ \frac{|\tilde \mP|}{q^{2k+2}}\Big(\frac{\delta^2}{8}\Big)^{2^m}
\ \ge\ \frac{1}{4q^2H^2}\Big(\frac{\delta^2}{8}\Big)^{2^m+1}.
\]
Define 
\[X(h)=\mathbb{E}_{a,b\in G_{k+1}}e(R(a,b,h)),\]
then:
\[\beta_1=\mE_{h\in G_{n-k}^m}X(h).\]
Define
\[\tilde{G}:=\{h:X(h)\geq\beta_1/2\}.\]

Since \(0\le X(h)\le1\) and \(\mathbb E_hX(h)=\beta_1\), this gives

$$ \frac{|\tilde G|}{q^{m(n-k)}}\ge\frac{\beta_1}{2}. $$

For every $h$, monotonicity of Gowers norms gives
\begin{align*}
    &X(h)=\|a\rightarrow e\circ d^mP(ah_1,\ldots, ah_m)\|_{U^1(G_{k+1})}^2\\
    &\le\|a\rightarrow e\circ d^mP(ah_1,\ldots, ah_m)\|_{U^m(G_{k+1})}^2
=\Big(\mathbb{E}_{\vec a\in G_{k+1}^m}e\big(\tilde Q(\vec a,h)\big)\Big)^{2^{1-m}},
\
\end{align*}

where
\[
\tilde Q(\vec a,h)=\sum_{\sigma\in S_m}d^mP\big(a_{\sigma(1)}h_1,\dots,a_{\sigma(m)}h_m\big)
\]
is $2m$-$\mF_p-$linear on $G_{k+1}^m\times G_{n-k}^m$. Hence
\[
\mathbb{E}_{\vec a}\,e(\tilde Q(\vec a,h))\ \ge\ \Big(\frac{\beta_1}{2}\Big)^{2^{m-1}}
\qquad\text{for } h\in\tilde G .
\]
For every other $h$, the map $\vec a\mapsto\tilde Q(\vec a,h)$ is $\mF_p-$multilinear  on the vector space
$G_{k+1}^m$. Therefore
\[
\mathbb{E}_{\vec a}e(\tilde Q(\vec a,h))
=\Pr_{a_2,\dots,a_m}\big[\tilde Q(\cdot,a_2,\dots,a_m,h)\equiv0\big]\ \ge0 .
\]
Keeping only $h\in\tilde G$,
\[
\mathbb{E}_{\vec a\in G_{k+1}^m,\ h\in G_{n-k}^m}e(\tilde Q(\vec a,h))
\ \ge\ \frac{|\tilde G|}{q^{m(n-k)}}\Big(\frac{\beta_1}{2}\Big)^{2^{m-1}}
\ \ge\ \Big(\frac{\beta_1}{2}\Big)^{2^{m-1}+1}=:\beta_2 .
\]
Observe that: 
\[\log_p\beta_2^{-1}=(2^{m-1}+1)\log_p(2/\beta_1)\le C(\delta,m,q,r)\log H\le C(r,\delta,m,q)\log\Lambda.\]
By Theorem \ref{BR} (over $\mF_p$), applied to the $2m$-linear form $\tilde Q$,
\[
r_{pr}(\tilde Q)\ \le\ C_{2m}\big(\log_p\beta_2^{-1}\big)^{C_{2m}}\ \le\ C_1\,(\log\Lambda)^{C_{2m}}=:s .
\]

 Write:
\[
\tilde Q=\sum_{i=1}^{s}X_i(x_{I_i})\,Y_i(x_{I_i^c}),
\]
where $x=(a_1,\dots,a_m,h_1,\dots,h_m)$, each $I_i$ is a nonempty proper subset of the $2m$ variables,
and $X_i,Y_i$ are multilinear. Restrict to the diagonal $a_1=\dots=a_m=a\in G_{k+1}\simeq\mF_p^{r(k+1)}$. Since $d^mP$ is
symmetric,
\[
\tilde Q(a,\dots,a,h)=m!\,d^mP(ah_1,\dots,ah_m).
\]
For each $i$ there are two cases.
\begin{itemize}
\item[(i)] One of $I_i,I_i^c$, say $I_i$, contains no $h$-variable. Then on the diagonal the term
equals $N_i(a)\,M_i(a,h)$, where $N_i(a)=X_i(a,\dots,a)$ depends only on $a$.
\item[(ii)] Both $I_i$ and $I_i^c$ contain $h$-variables. Then for every fixed $a$ the term is a
product of two multilinear forms in disjoint nonempty sets of $h$-variables, so it has partition rank
$\le1$ in $h$.
\end{itemize}
Each $N_i$ from case~(i) is a polynomial of degree $\le m$ in the $r(k+1)$ $\mathbb{F}_p$-coordinates
of $a\in G_{k+1}$. There are at most $s$ of them, and $a=0$ is a common zero.

For large $n$ we have $r(k+1)>rW=r\Lambda^{1/8}>rms=rmC_1(\log\Lambda)^{C_{2m}}$,
hence we can use Chevalley–Warning to find a common zero $a\ne0$. For this $a$, only the case~(ii) terms survive, and since
$m!\in\mathbb{F}_p^\times$:
\[
r_{pr}\big(h\mapsto d^mP(ah_1,\dots,ah_m)\big)\ \le\ s .
\]
Since $a\neq 0$ and $\deg(a)\leq k$, we can use Lemma \ref{sub} and deduce:
\[r_{pr}(d^mTr(t\tilde{P}))\leq rmk+s\leq C_1(\log \Lambda)^{C_{2m}}+rmH\]
\end{proof}

Fix $\delta\in(0,1)$ and let $n\ge n_1$. 

Suppose that
\[|\mathbb{E}_{g\in G_n}\nu(g)\alpha(\tilde{P}(g))|>\delta\] for some $\tilde{P}$ of degree $\le m$ and some nontrivial
$\alpha$. By Lemma~\ref{key} and Remark~\ref{r-ineqaulity} (polarisation, $p>m$),
\[
Tr(\tilde{tP})=\sum_{i\le s'}S_iT_i+L,\qquad s'\le C_1(\log\Lambda)^{C_{2m}}+rmH,
\]
with all $S_i,T_i,L$ of degree $\le m-1$. To complete the proof, we need the following Lemma: 
\begin{Lemma }\label{FA}.
Assume that for every polynomial $P$ of degree $\leq m-1$:
\[|\mE_{g\in G_n}\nu(g)e_p(-P(g))|\leq q^{-n}\phi(n).\]
Then for every $s\in \mN$, every $1$-bounded function $F:\mF_q^s \xrightarrow{} \mC$, and every $P_{1},\ldots ,P_{s}$ of degree $\leq m-1$:
\[|\mE_{g\in G_n}\nu(g)F(P_{1}(g),\ldots,P_{s}(g))|\leq q^{s-n}\phi(n).\]
\end{Lemma }
Lemma \ref{FA} is an easy corollary of a similar Lemma in \cite{mobius}. We give the proof at the end of the paper, for the convenience of the reader. \\

Given $\alpha=e_p(t\cdot )$, apply Lemma \ref{FA} with $2s'+1$ arguments and
$F(x,y,z)=e\big(\sum_ix_iy_i+z\big)$ on $\mF_p^{2s'+1}$ and extended to zero to $\mF_q^{2s'+1}$:
\begin{align*}
    &\Big|\mathbb{E}_{g\in G_n}\nu(g)\alpha(\tilde{P}(g))\Big|=\Big|\mathbb{E}_{g\in G_n}\nu(g)e(P(g))\Big|\\
    &\le q^{2s'+1}\frac{\phi_{m-1}(n)}{q^n}
\le q^{\,2C_1(\log\Lambda)^{C_{2m}}+2mr\Lambda^{1/4}+1-\Lambda}.
\end{align*}

The exponent tends to $-\infty$, depending only on $(\delta,m,q,\Lambda(n))$. So for $n\ge n_0(\delta)$
the right-hand side is $<\delta$, a contradiction. Hence
$\sup_{\alpha}\sup_{\deg P\le m}|\mathbb{E}_{G_n}\nu\,\alpha(P)|\le\delta$ for all $n\ge n_0(\delta)$. As $\delta$ was arbitrary, $\phi_m(n)=o_{m,\nu}(q^n)$.

\section{Halasz Criterion for pretentious distance}
We are now interested in characterizing the class of aperiodic multiplicative functions in term of pretentious-distance (definition \ref{pd}).


In the context of $\mathbb{N}$, a similar characterization exists (see, for example, \cite{HOFA}), where aperiodic functions are exactly those that do not ``pretend'' to be twisted Dirichlet characters; i.e., $\min_{\chi, \theta} \mathbb{D}(f, \chi e_{\theta}, \infty) = \infty$. We aim to provide an analogous characterization in the function field setting.





To bridge the mean-value results with the pretentious framework, we invoke the following version of the Hal\'{a}sz--Delange Theorem for function fields, as established in \cite{Klurman thesis}, Theorem 1.4.1:
\begin{theorem}[Hal\'{a}sz--Delange Theorem for function fields]\label{haldel}
    For a $1$-bounded multiplicative function $f:\mF_q[x]\rightarrow\mC$, the following holds:
    \begin{enumerate}
        \item Let $N$. Define $\mD_f(N)=\min_{\theta}\mD(f,e_{\theta},N).$ Then:
        \[\mE_{g \text{ monic},\deg(g)=N}f(g)\ll (1+\mD_f(N))e^{-\mD_f(N)}\]
        \item There exists some $\theta_0\in[0,1)$ such that $\mD(f,e_{\theta_0},\infty)<\infty$. For $\epsilon>0$ denote $m(n)=\lceil(1-\epsilon)\frac{\log(n)}{\log(q)}\rceil$. Then:
        \[\mE_{g \text{ monic}, \deg(g)=n}f(g)=e(2\pi i \theta_0n)P(fe_{-\theta_0},n)+O_{\epsilon}\!\left(\mD(fe_{\theta_0},m(n),n)+\frac{1}{n^{1-\epsilon}}\right),\]
        where:
        \begin{align*}
            &P(f,n)=\prod_{r \text{ irred., } \deg(r)\leq n}\left(1-q^{-\deg(r)}\right)\!\left(\sum_{k=0}^{\infty}f(r^k)q^{-k\deg(r)}\right) \\
            &\text{ and }\mD(fe_{\theta_0},m(n),n)^2=\mD(f,e_{\theta_0},n)^2-\mD(f,e_{\theta_0},m(n))^2.
        \end{align*}
    \end{enumerate}
\end{theorem}
Using Theorem \ref{haldel}, we establish the link between aperiodicity and pretentiousness. 
\begin{theorem}\label{apf-equiv-np}
    Let $\nu$ a multiplicative function bounded by $1$. Then $\nu$ is aperiodic if and only if $\nu$ is non-pretentious.
\end{theorem}
\begin{proof}
\textit{ Aperiodicity$\impliedby$ non-pretentious.} Let $\nu$ be a non-pretentious function, let $Q\in\mF_q[x]$ and $L\in\mN$, let $n$ and let and let $a$ be a $Q,L$-Hayes-periodic function, $\deg(Q)=K$, i.e for every $g,\deg(g)=n$ $a(g)=a(g\text{ mod Q}, g_{n-L},\ldots,g_n)$. We aim to prove that: 
\[\sum_{\deg g=n}\nu(g)a(g)=o(q^n).\]

Let $\mathcal{M}_n$ be the set of polynomials of degree exactly $n$. Let $t=x^{-1}$ and:
\[
U_L=\big(\mathbb{F}_q[t]\big)\big/\big(t^{L+1}\big),\qquad
\rho_L(g)=g\,x^{-\deg g}\bmod x^{-L-1}.
\]
Then $U_L$ is a ring of order $q^{L+1}$, the map $\rho_L$ is multiplicative.\\

For $g\in\mathcal{M}_n$ with $n\geq L$, $\rho_L(g)$  is determined by
$(g_{n},\dots,g_{n-L})$. A Hayes character modulo $(Q,L)$ is a character $\tilde\chi$ of
$\Gamma_{Q,L}=(\mathbb{F}_q[x]/Q)^\times\times U_L^\times$, evaluated at $(g\bmod Q,\rho_L(g))$ and set equal to
$0$ if $(g,Q)\neq1$; its conductor is at most $K+L$.\\

Define:
\[C_n(r,c)=\{g,\deg(g)=n, g\equiv_Q r, \rho_L(g)=c\}.\]
Using this definition, it's enough to prove that for every $(r,c)\in U_L$:
\[\sum_{g\in C_n(r,c)}\nu(g)=o(q^n).\]

For each $r$ with $(r,Q)=1$, the function $1_{g \in C_n{(r,c)}}$ is a finite combination of the characters of $\Gamma_{Q,L}$, and hence a finite combination of Hayes-character with conductor $K+L$. Meaning, it's enough to prove that $\nu$ is orthogonal to each Hayes-character $\tilde{\chi}$ of this form. Since:
\begin{align*}
    &\min_{\theta}\mD(\nu\tilde{\chi},e_{\theta},N)=\min_{\theta}\mD(\nu,\tilde{\chi}e_{\theta},N)\\
    &\geq \min_{\tilde{\chi},cond(\tilde{\chi})\leq K+L}\min_{\theta}\mD(\nu,\tilde{\chi}e_{\theta},N)\rightarrow \infty
\end{align*}
Since:
\[\sum_{g,\deg(g)=n}\nu(g)=\sum_{u\in\mF_q^{\times}}\nu(u)\sum_{g \text{ monic, }\deg(g)=n}\nu(g)\]
Using \ref{haldel}, we deduce
\begin{align*}
    &\sum_{g,\deg(g)=n}\nu(g)\tilde{\chi}(g)=o(q^n),
\end{align*}
And hence:
\begin{align*}
    &\sum_{g,\deg(g)=n}\nu(g)1_{g \in C_n{(r,c)}}=o(q^n).
\end{align*}

If $(r,Q)=b\neq 1$, then $(r/b,Q/b)=1$. Factor every $g\in\mathcal{M}_n$ uniquely as $g=uv$, where
$u$ is monic, every prime factor of $u$ divides $b$, and $(v,b)=1$. Since $(u,v)=1$, we have
$\nu(g)=\nu(u)\nu(v)$. Let
\[
A_n=\{u\ \text{monic}:\ \deg(u)\le n,\ \text{every prime factor of $u$ divides $b$}\}.
\]

\emph{The congruence condition.} We have $g\equiv_Q r$ if and only if $b\mid u$ and
$v\frac ub\equiv_{Q/b}\frac rb$. Indeed, if $g\equiv_Q r$ then $b\mid g$, hence $b\mid u$ because $(v,b)=1$,
and dividing by $b$ gives the second condition.

Conversely, multiply by $b$. Since $(\frac rb,\frac Qb)=1$,
the condition $v\frac ub\equiv_{Q/b}\frac rb$ forces $(v,\frac Qb)=1$, hence $(v,Q)=1$, and forces
$(\frac ub,\frac Qb)=1$; if $(\frac ub,\frac Qb)\neq1$ there is no such $v$. So, for fixed $u$, the
admissible residues of $v$ form a union of at most $\Phi(Q)$ reduced residue classes modulo $Q$.

\emph{The top-coefficient condition.} Since $\rho_L$ is multiplicative and $U_L$ is a ring,
\[
\rho_L(uv)=c\quad\Longleftrightarrow\quad \rho_\ell(v)=c\,\rho_L(u)^{-1}.
\]

Consequently, for fixed $u\in A_n$ with $b\mid u$, the set of admissible $v\in\mathcal{M}_{n-\deg u}$ is a
disjoint union of at most $\Phi(Q)$ sets $C_{n-\deg u}(r',c')$ with $(r',Q)=1$ and
$c'=c\,\rho_\ell(u)^{-1}$. For such $u$ put
\[
\varepsilon_u(m)=\frac1{q^{m}}\sum_{\substack{v\in\mathcal{M}_m,\ (v,b)=1\\
v\frac ub\equiv_{Q/b}\frac rb,\ \ \rho_\ell(v)=c\,\rho_\ell(u)^{-1}}}\nu(v).
\]
Then $|\varepsilon_u(m)|\le1$, and by the co-prime case, $\varepsilon_u(m)\to0$ as $m\to\infty$ for each fixed $u$. By the
discussion above,
\[
\frac1{q^n}\sum_{g\in C_n(r,c)}\nu(g)=\sum_{\substack{u\in A_n\\ b\mid u}}\nu(u)\,q^{-\deg u}\,\varepsilon_u(n-\deg u).
\]
Let $s$ be the number of distinct prime factors of $b$; there are at most $(k+1)^s$ elements of $A_n$ of
degree $k$. Fix $T\ge\deg b$. The terms with $\deg u>T$ contribute at most $\sum_{k>T}(k+1)^s/^{-k}$. The
terms with $\deg u\le T$ are finitely many and do not depend on $n$; for each of them
$n-\deg u\to\infty$, so they tend to $0$. Hence
\[
\limsup_{n\to\infty}\Big|\frac1{q^n}\sum_{g\in C_n(r,c)}\nu(g)\Big|\le\sum_{k>T}(k+1)^sq^{-k},
\]
and letting $T\to\infty$ gives the desire result.


\textit{Aperiodicity $\implies$ non-pretentious.} Assume $\nu$ is aperiodic. 

Observe that by definition, $\nu$ is orthogonal to every Hayes-periodic function, and in
particular to every Hayes character $\tilde\chi=\chi\xi$.

 We show that
$\mD(\nu,\tilde\chi e_\theta,\infty)=\infty$ for every Hayes character $\tilde\chi$ and every
$\theta\in[0,1)$.

Suppose, for contradiction, that $\mD(\nu,\tilde\chi_0e_{\theta_0},\infty)<\infty$, where $\tilde\chi_0$ is
a Hayes character modulo $(b,\ell)$.  Put
\[
\tilde\nu(g)=\nu(g)\,\overline{\tilde\chi_0(g)}\,1_{(g,Q_0)=1},\qquad
Q_0=\begin{cases}x(x+1),&q=2,\\ 1,&q\ge3.\end{cases}
\]
Then $\tilde\nu$ is $1$-bounded and multiplicative, and $\mD(\tilde\nu,e_{\theta_0},\infty)<\infty$,
since the factor $1_{(g,Q_0)=1}$ changes the distance at no more than two primes.

\emph{(i) The mean of $\tilde\nu$ tends to $0$.} The function $\overline{\tilde\chi_0}\,1_{(\cdot,Q_0)=1}$
is Hayes-periodic (modulus $bQ_0$, top coefficients of length $\ell$). Since $\nu$ is aperiodic,
\[
\mE_{g,\deg(g)=n}\tilde\nu(g)=q^{-n}\sum_{0\neq r\in\mF_q}\nu(r)\chi_0(r)\sum_{g \text{ monic },deg(g)=n}\tilde\nu(g)\xrightarrow[n\to\infty]{}0.
\]

\emph{(ii) Hal\'asz--Delange.} By Theorem~\ref{haldel}(2) applied to $\tilde\nu$ with $\epsilon=\frac12$,
\[
\mE_{g,\deg(g)=n}\tilde\nu(g)=e(\theta_0n)\,P(\tilde\nu e_{-\theta_0},n)
+O\!\left(\mD(\tilde\nu,e_{\theta_0},m(n),n)+n^{-1/2}\right).
\]
Since $\mD(\tilde\nu,e_{\theta_0},\infty)<\infty$ and $m(n)\to\infty$, the tail
$\mD(\tilde\nu,e_{\theta_0},m(n),n)$ tends to $0$. Together with (i), this gives
\[
|P(\tilde\nu e_{-\theta_0},n)|\xrightarrow[n\to\infty]{}0 .
\]

\emph{(iii) The Euler product is bounded below.} Write $\nu_0=\tilde\nu e_{-\theta_0}$, and for an irreducible
$r$ let
\[
L_r=\big(1-q^{-\deg(r)}\big)\sum_{k\ge0}\frac{\nu_0(r^k)}{q^{k\deg(r)}}
\]
be the local factor of $P(\nu_0,n)$ at $r$. If $q^{\deg(r)}\ge8$, then
\[
L_r=1+z_r,\qquad z_r=\frac{\nu_0(r)-1}{q^{\deg(r)}}+O\big(q^{-2\deg(r)}\big),\qquad |z_r|\le\tfrac12,
\]
and therefore $\log|L_r|\ge-\frac{1-\operatorname{Re}\nu_0(r)}{q^{\deg(r)}}-\frac{C}{q^{2\deg(r)}}$.
For the finitely many $r$ with $q^{\deg(r)}<8$:
\begin{itemize}
\item if $q^{\deg(r)}\ge3$, then $|L_r|\ge\big(1-q^{-\deg(r)}\big)\big(1-\frac1{q^{\deg(r)}-1}\big)>0$;
\item if $q^{\deg(r)}=2$, then $r\mid Q_0$, so $\nu_0(r^k)=0$ for $k\ge1$ and $L_r=\frac12$.
\end{itemize}
Since $1-\operatorname{Re}\nu_0(r)=1-\operatorname{Re}\big(\tilde\nu(r)e_{-\theta_0}(r)\big)$, summing over
$r$ gives, for all $n$,
\[
\log|P(\nu_0,n)|\ \ge\ -\,\mD(\tilde\nu,e_{\theta_0},\infty)^2-C'\ >\ -\infty .
\]
This contradicts (ii). Hence $\mD(\nu,\tilde\chi e_\theta,\infty)=\infty$ for every $\tilde\chi$ and
every $\theta$.

\begin{remark}\label{theta-uniform}
In this case:
\[
\min_{\theta\in[0,1]}\mD(\nu,e_\theta,N)\xrightarrow[N\to\infty]{}\infty .
\]
\begin{proof}
  Assume by contradiction we have a sequence $\theta_j,N_j$ with $N_j\rightarrow\infty$ and $M$ such that:
  \[\mD(\nu,\theta_j,N_j)\leq M.\]
Take a subsequence of $\theta_j$ that is convergence to $\theta_0$. Fix some $N_0$, then by monotonicity in $N$ we have that for large enough $j$:
\[\mD(\nu,\theta_j,N_0)\leq \mD(\nu,\theta_j,N_j)\leq M.\]
Since $N_0$ is fixed, the map $\theta\rightarrow\mD(\nu,\theta,N_0)$ is continuous. Hence:
  \[\mD(\nu,\theta_0,N_0)=\lim_k\mD(\nu,\theta_{j_k},N)\leq M.\]
  This equation holds for every $N_0$, hence it contradicts point-wise convergence at $\theta_0$.
\end{proof}
\end{remark}

Finally, Fix $K$, and let $\tilde\chi$ be a Hayes character with
$\operatorname{cond}(\tilde\chi)\le K$. By Step~1, $\mD(\nu\overline{\tilde\chi},e_\theta,\infty)=\infty$
for every $\theta$. Using ~\ref{theta-uniform}, applied to $\nu\overline{\tilde\chi}$ in place of $\nu$,
gives $\min_\theta\mD(\nu,\tilde\chi e_\theta,N)\to\infty$. There are only finitely many Hayes characters
of conductor $\le K$, so
\[
\lim_{N\to\infty}\ \min_{\operatorname{cond}(\tilde\chi)\le K}\ \min_{\theta}\ \mD(\nu,\tilde\chi e_\theta,N)=\infty .
\]
That is, $\nu$ is non-pretentious. 
\end{proof}

\section{Appendix}\label{appendix}
\subsection{Katai's criterion}
In this section, we prove Katai's criterion for $\mathbb{F}_q[x]$. The proof follows the same strategy as the analogue in $\mathbb{N}$ (see for example \cite{TTB}).\\

Throughout, $r,s$ denote \emph{monic irreducible} polynomials, $P_k$ is the set of monic
irreducibles of degree $k$, and $\pi(k)=|P_k|$. We use the standard bounds
\begin{equation}\label{eq:PPT}
\frac{q^k}{2k}\le \pi(k)\le \frac{q^k}{k}\qquad(k\ge 4).
\end{equation}

For integers $0\le W<H$ put:
\[
A=A(W,H)=\sum_{W<\deg(r)<H}\frac{1}{q^{\deg(r)}}=\sum_{W<k<H}\frac{\pi(k)}{q^k}.
\]
By \eqref{eq:PPT}, if $W\ge3$ then $A\ge\frac12\sum_{W<k<H}\frac1k\ge\frac12\log\frac{H}{W+1}$.

\begin{Lemma }[Tur\'an--Kubilius for function fields]\label{TK}
Let $0\le W<H$ be integers with $2H\le n$. Then
\begin{equation}\label{k}
\sum_{g\in G_n}\Big|\sum_{\substack{W<\deg(r)<H\\ r\mid g}}1-A\Big|^2\le A\,q^n.
\end{equation}
\end{Lemma }

\begin{proof}
Write $\omega(g)=\sum_{W<\deg(r)<H,\ r\mid g}1$. For $\deg(r)\le n$, the number of $g\in G_n$
divisible by $r$ is exactly $q^{n-\deg(r)}$. Hence
\[
\sum_{g\in G_n}\omega(g)=\sum_{W<\deg(r)<H}q^{n-\deg(r)}=q^nA.
\]
Similarly,
\[
\sum_{g\in G_n}\omega(g)^2=\sum_{W<\deg(r),\deg(s)<H}\ \sum_{\substack{g\in G_n\\ r,s\mid g}}1,
\qquad
\sum_{\substack{g\in G_n\\ r,s\mid g}}1=\begin{cases}q^{n-\deg(r)}& r=s,\\
q^{n-\deg(r)-\deg(s)}& r\neq s,\end{cases}
\]
where for $r\neq s$ we used that $r,s$ are coprime and $\deg(rs)<2H\le n$. Therefore
\[
\sum_{g\in G_n}\omega(g)^2\le q^nA+q^nA^2,
\]
and
\[
\sum_{g\in G_n}\big(\omega(g)-A\big)^2=\sum_g\omega^2-2A\sum_g\omega+A^2q^n
\le q^nA+q^nA^2-2q^nA^2+q^nA^2=q^nA.\qedhere
\]
\end{proof}

\begin{theorem}[K\'atai's criterion, quantitative form]\ref{p_katai}
Let $f:\mF_q[x]\to\mathbb{C}$ with $|f|\le1$, and let $\nu:\mF_q[x]\to\mathbb{C}$ be
multiplicative with $|\nu|\le1$. Let $0\le W<H$ be integers with $2H\le n$, and let $\eta\ge0$
be such that for every integer $k$ with $W<k<H$,
\begin{equation}\label{eq.4}
\sum_{a,b\in P_k}\Big|\sum_{g\in G_{n-k}}f(ag)\overline{f}(bg)\Big|\le \eta\,|P_k|^2q^{n-k}.
\end{equation}
Then
\[
\Big|\sum_{g\in G_n}\nu(g)f(g)\Big|\le q^n\Big(A^{-1/2}+\frac{2q^{-W}}{(W+1)A}+\eta^{1/2}\Big).
\]
\end{theorem}

\begin{proof}
Write $\omega(g)=\sum_{W<\deg(r)<H,\ r\mid g}1$.

 By Cauchy--Schwarz and Lemma~\ref{TK},
\[
\Big|\sum_{g\in G_n}\big(\omega(g)-A\big)\nu(g)f(g)\Big|\le (q^n)^{1/2}(Aq^n)^{1/2}=A^{1/2}q^n.
\]
Rearranging,
\begin{equation}\label{eq:step1}
\Big|\sum_{g\in G_n}\nu(g)f(g)\Big|\le\frac1A\Big|\sum_{g\in G_n}\ \sum_{W<\deg(r)<H}\nu(g)f(g)1_{r\mid g}\Big|+A^{-1/2}q^n.
\end{equation}

 Every $g\in G_n$ with $r\mid g$ is $g=rh$ for a unique $h\in G_{n-\deg(r)}$.
If $r\nmid h$ then $\nu(g)f(g)=\nu(r)\nu(h)f(rh)$. The remaining $g$ satisfy $r^2\mid g$; there
are $q^{n-2\deg(r)}$ of them (as $n-\deg(r)\ge\deg(r)$), and for each
$|\nu(g)f(g)-\nu(r)\nu(h)f(rh)|\le2$. By \eqref{eq:PPT}, the total contribution of these
exceptions is at most
\[
2\sum_{W<k<H}\pi(k)q^{n-2k}\le 2q^n\sum_{k\ge W+1}\frac{q^{-k}}{k}\le\frac{2q^{n-W}}{W+1}.
\]
Hence
\begin{equation}\label{eq:step2}
\Big|\sum_{g}\sum_{W<\deg(r)<H}\nu(g)f(g)1_{r\mid g}\Big|
\le\sum_{W<k<H}\Big|\sum_{r\in P_k}\nu(r)\sum_{h\in G_{n-k}}\nu(h)f(rh)\Big|+\frac{2q^{n-W}}{W+1}.
\end{equation}

 Fix $k$ with $W<k<H$. By Cauchy--Schwarz in $h$,
\[
\Big|\sum_{h\in G_{n-k}}\nu(h)\sum_{r\in P_k}\nu(r)f(rh)\Big|^2
\le q^{n-k}\sum_{h\in G_{n-k}}\Big|\sum_{r\in P_k}\nu(r)f(rh)\Big|^2,
\]
and
\begin{align*}
    &\sum_{h\in G_{n-k}}\Big|\sum_{r\in P_k}\nu(r)f(rh)\Big|^2
    =\sum_{r,s\in P_k}\nu(r)\overline{\nu}(s)\sum_{h\in G_{n-k}}f(rh)\overline{f}(sh)\\
    &\leq \sum_{r,s\in P_k}\Big|\sum_{h\in G_{n-k}}f(rh)\overline{f}(sh)\Big|
\le\eta\,|P_k|^2q^{n-k}
\end{align*}

by \eqref{eq.4}. Hence
\[
\Big|\sum_{r\in P_k}\nu(r)\sum_{h\in G_{n-k}}\nu(h)f(rh)\Big|\le\eta^{1/2}|P_k|\,q^{n-k},
\]
and summing over $k$,
\[
\sum_{W<k<H}\Big|\sum_{r\in P_k}\nu(r)\sum_{h}\nu(h)f(rh)\Big|\le\eta^{1/2}q^n\sum_{W<k<H}\frac{\pi(k)}{q^k}=\eta^{1/2}A\,q^n.
\]

Inserting into \eqref{eq:step2} and then into \eqref{eq:step1},
\[
\Big|\sum_{g\in G_n}\nu(g)f(g)\Big|\le\frac1A\Big(\eta^{1/2}Aq^n+\frac{2q^{n-W}}{W+1}\Big)+A^{-1/2}q^n
= q^n\Big(A^{-1/2}+\frac{2q^{-W}}{(W+1)A}+\eta^{1/2}\Big).\qedhere
\]
\end{proof}

\begin{remark}
The diagonal terms $a=b$ in \eqref{eq:katai-hyp} contribute $|P_k|q^{n-k}$, so necessarily
$\eta\ge\max_{W<k<H}|P_k|^{-1}\ge q^{-W}$. For the bound to tend to $0$ one therefore needs
$W\to\infty$, as well as $A\to\infty$.
\end{remark}

\begin{corollary}[K\'atai's criterion, qualitative form]
Let $\phi:\mathbb{N}\to\mathbb{N}$ with $\phi(n)\to\infty$ and $\phi(n)+\phi(n)^2\le n/2$, and set
$W=\phi(n)$, $H=\phi(n)+\phi(n)^2$. Assume
\[
\eta(n):=\max_{W<k<H}\ \frac{1}{|P_k|^2q^{n-k}}\sum_{a,b\in P_k}\Big|\sum_{g\in G_{n-k}}f(ag)\overline{f}(bg)\Big|\ \xrightarrow[n\to\infty]{}\ 0.
\]
Then $\sum_{g\in G_n}\nu(g)f(g)=o(q^n)$ for every multiplicative $\nu$ with $|\nu|\le1$.
\end{corollary}

\begin{proof}
Here $A\ge\frac12\log\frac{\phi+\phi^2}{\phi+1}\to\infty$ (for $\phi\ge3$). Apply
Theorem~\ref{p_katai}: all three terms tend to $0$.
\end{proof}

\subsection{Proof of Lemma \ref{FA}}
In this section, we prove the following lemma:
\begin{Lemma *}[\ref{FA}, restated]
Assume that for every polynomial $P$ of degree $\leq m-1$:
\[|\mE_{g\in G_n}\nu(g)e_p(-P(g))|\leq q^{-n}\phi(n).\]
Then for every $s\in \mN$, every $1$-bounded function $F:\mF_q^s \xrightarrow{} \mC$, and every $P_{1},\ldots ,P_{s}$ of degree $\leq m-1$:
\[|\mE_{g\in G_n}\nu(g)F(P_{1}(g),\ldots,P_{s}(g))|\leq q^{s-n}\phi(n).\]
\end{Lemma *}
Its proof is similar to the proof of the analogous lemma in \cite{mobius}.
\begin{proof}
For every $a=(a_1,\ldots,a_s)\in \mF_{q}^s$, denote 
\[
V_a=\{g\in G_n\mid P_{i}(g)=a_i \text{ for every } i\}.
\]
Then:
\[
\sum_{g\in G_n}\nu(g)F(P_{1}(g),\ldots,P_{s}(g))=\sum_{a\in \mF_{q}^s}F(a_1,\ldots,a_s)\sum_{g\in V_a}\nu(g).
\]
Note that 
\[
1_{V_a}(g)=\mE_{\chi=(\chi_1,\ldots,\chi_s)\in \widehat{\mF_{q}^s}}\prod_{i\in[s]} \chi_i(P_i(g)-a_i),
\]
hence 
\[
\sum_{g\in V_{a}}\nu(g)=\mE_{\chi\in \widehat{\mF_{q}^s}}\prod_{i \in [s]}\chi_i(-a_i)\sum_{g \in G_n}\nu(g)\prod_{i\in [s]}\chi_i(P_i(g)).
\]
By the triangle inequality, 
\[
\left|\sum_{g\in V_{a}}\nu(g)\right|\leq \max_{\chi}\left|\sum_{g\in G_n}\nu(g)\prod_{i \in [s]}\chi_i(P_i(g))\right|.
\] 
Since each $\chi_i(x)=e_p(t_i x)$ for some $t_i$:
\[
\prod_{i\in [s]}\chi_i(P_i(g))=e_p\!\left(\sum_{i\in[s]} t_i P_i(g)\right),
\] 
where $\sum_{i\in[s]} t_i P_i(g)$ is a polynomial of degree $\leq m-1$.
Hence:
\[
\left|\sum_{g\in V_{a}}\nu(g)\right|\leq \phi(n).
\]
It follows that:
\begin{align*}
    \left|\mE_{g\in G_n}\nu(g)F(P_{1}(g),\ldots,P_{s}(g))\right|
    &=q^{-n}\left|\sum_{a\in \mF_{q}^s}F(a_1,\ldots,a_r)\sum_{g \in V_a}\nu(g)\right|
    \leq q^{s-n}\phi(n),
\end{align*}
as desired.
\end{proof}

\end{document}